\documentclass[a4paper,twoside]{article}
\usepackage{amsmath,amsfonts,amssymb,amsthm}
\usepackage{enumerate} 
\usepackage[bookmarks]{hyperref}
\usepackage[center]{titlesec}
\usepackage{
mathrsfs,dsfont}
\usepackage{tocbibind}
\usepackage{geometry}
\numberwithin{equation}{section}

\newtheorem{thm}{Theorem}[section]
\newtheorem{cj}{Conjecture}[section]
\newtheorem{defi}[thm]{Definition}
\newtheorem{pro}[thm]{Proposition}
\newtheorem{lem}[thm]{Lemma}
\newtheorem{cor}[thm]{Corollary}
\newtheorem{rem}[thm]{Remark}

\def\bthm{\begin{thm}\def\ethm{\end{thm}}}
\def\bpro{\begin{pro}\def\epro{\end{pro}}}

\def\blem{\begin{lem}\def\elem{\end{lem}}}

\def\bitm{\begin{itemize}} \def\eitm{\end{itemize}}
\def\benu{\begin{enumerate}} \def\eenu{\end{enumerate}}

\def\bpf{\begin{proof}}\def\epf{\end{proof}}
\def\beq{\begin{equation}}\def\eeq{\end{equation}}
\def\beqs{\begin{eqnarray}}\def\eeqs{\end{eqnarray}}
\def\beqsnl{\begin{eqnarray*}}\def\eeqsnl{\end{eqnarray*}}

\def\bb{\mathbb}
\def\ca{\mathcal}
\def\fr{\mathfrak}

\def\Xint#1{\mathchoice
  {\XXint\displaystyle\textstyle{#1}}%
  {\XXint\textstyle\scriptstyle{#1}}%
  {\XXint\scriptstyle\scriptscriptstyle{#1}}%
  {\XXint\scriptscriptstyle\scriptscriptstyle{#1}}%
  \!\int}
\def\XXint#1#2#3{{\setbox0=\hbox{$#1{#2#3}{\int}$}
  \vcenter{\hbox{$#2#3$}}\kern-.5\wd0}}

\def\dashint{\Xint-}

\usepackage{fancyhdr}
\newcommand{\bra}[1]{\left(#1\right)}

\newcommand{\abs}[1]{\left\vert#1\right\vert}
\newcommand{\set}[1]{\left\{#1\right\}}
\newcommand{\ip}[1]{\left<#1\right>}
\newcommand{\norm}[1]{\left\Vert#1\right\Vert}

\DeclareMathOperator{\re}{Re}
\DeclareMathOperator{\im}{Im}
\DeclareMathOperator{\s}{\bb{S}^{2n+1}}
\DeclareMathOperator{\h}{\bb{H}^n}
\DeclareMathOperator{\m}{\ca{M}}

\begin{document}
\title{\textbf{Remainder Terms for 
Several Inequalities on Some Groups of Heisenberg-type
}}
\footnotetext{\emph{Date:} March. 28, 2014, revised on Aug. 29, 2014.}
\footnotetext{
Heping Liu is supported by National Natural Science Foundation of China under Grant No. 11371036 and the Specialized Research Fund for the Doctoral Program of Higher Education of China under Grant No. 2012000110059. An Zhang is supported by China Scholarship Council under Grant No. 201306010009.}

\author{Heping Liu, An Zhang\footnote{Corresponding author.}}

\date{}
\maketitle
\begin{abstract}
We give estimates of the remainder terms for several conformally-invariant Sobolev-type inequalities on the Heisenberg group, in analogy with the Euclidean case. By considering the variation of associated functionals, we give a \emph{stability} of two dual forms: the fractional Sobolev (Folland-Stein) and Hardy-Littlewood-Sobolev inequality, in terms of distance to the submanifold of \emph{extremizers}. Then we compare their remainder terms to improve the inequalities in another way. We also compare, in the limit case $s=Q$ (or $\lambda=0$), the remainder terms of Beckner-Onofri inequality and its dual Logarithmic Hardy-Littlewood-Sobolev inequality. Besides, we also list without proof some results for the other two cases of groups of Iwasawa-type. Our results generalize earlier works on Euclidean spaces by Chen, Frank, Weth \cite{cfw13} and Dolbeault, Jankowiakin \cite{dj14} onto some groups of Heisenberg-type.
\end{abstract}

\section{Introduction}
In this paper, we consider on Iwasawa-type groups (for which, the Heisenberg group $\bb{H}^n$ is the simplest non-Euclidean one) the problem of sharpening the conformally-invariant Sobolev-type inequalities obtained in \cite{fl12,clz14,clz14-}, by adding a remainder term proportional to distance square to the submainifold of extremizers, motivated originally from the question asked by Brezis and Lieb in \cite{bl85} for Euclidean spaces. We also want to compare the remainder terms of the dual inequalities similar to provious works on Euclidean spaces like that in \cite{dj14}.

On Euclidean case $\bb{R}^n$, we have classical fractional Sobolev inequality: given any exponent $0<s<n, q=\frac{2n}{n-s}$, for any function $f$ in $\frac{s}{2}$-order (homogeneous) Sobolev space $\dot{H}^{s/2}$ endowed with the norm $\norm{f}_{\dot{H}^{s/2}}=\norm{(-\Delta)^{s/4}f}_{L^2}=(\int f (-\Delta)^{s/2}f)^{1/2}$, we have
\beq\label{fs-e}\norm{f}_{\dot{H}^{s/2}}\gtrsim \norm{f}_{L^q},\eeq
and the sharp constant was first computed by several pioneers for some special cases and finally obtained by Lieb \cite{lie83} (there Lieb consider a dual form, the (diagonal) Hardy-Littlewood-Sobolev inequality) for exponent $s$ and dimension $n$ of all range and the author also proved that equality for sharp (\ref{fs-e}) can only be achieved by extremal functions (called \emph{extremizers}) of the form $c(1+|\delta x-x_0|^2)^{-\frac{n-s}{2}}, c\in\bb{R}\setminus \{0\},\delta>0, x_0\in\bb{R}^n$. Note that, the inequality (\ref{fs-e}) has an equivalent edition on the $n$-sphere $\bb{S}^n$ through the stereographic transformation $\ca{S}: x\in\bb{R}^n \mapsto \zeta=\bra{\frac{2x}{1+|x|^2},\frac{1-|x|^2}{1+|x|^2}} \in \bb{S}^n$. Both the fractional and its dual Hardy-Littlewood-Sobolev inequalities are invariant under one kind of ``conformal actions", including not only the translations, dilations, but also an inversion, which is given for the fractional inequality by $\sigma_{inv}: f(x)\mapsto |x|^{s-n}f(x|x|^{-2})$. The sharp inequalities and extremizers were got by exploiting this big conformal symmetry group in several beautiful ways. Then, in \cite{bl85}, Brezis and Lieb proposed a stability problem for the sharp inequality (\ref{fs-e}) in case $s=2$, asking whether there exists a positive constant $\alpha$, s.t., the following estimate holds:
\beq\label{sfs-e} \norm{f}_{\dot{H}^{s/2}}^2-C_{sharp} \norm{f}_{L^q}^2 \ge \alpha d^2(f,\ca{M}),\eeq
where $\ca{M}$ is the $(n+2)$-dimensional smooth submanifold of $\dot{H}^{s/2}$, consisting of all real-valued extremizers, and $d(f,\ca{M})$ is the usual distance of $f$ to $\ca{M}$ under the Sobolev norm. Bianchi and Egnell gave a positive answer in \cite{be91} still for $s=2$, which was later extended to more fractions but not of total range in \cite{lw00,bww03}.
Chen, Frank and Weth in \cite{cfw13} extended the method of Bianchi and Egnell to obtain above remainder term inequality (\ref{sfs-e}) for all fractions $0<s<n$, which contains all old results above. Naturally, we guess stability like (\ref{sfs-e}) should also hold for analogous inequalities (generalized fractional Sobolev, sometimes also called Folland-Stein, and Hardy-Littlewood-Sobolev inequality) on the Heisenberg and some more general groups of Heisenberg type.
In a different way, Dolbeault pointed out in \cite{dol11} that in $s=2$ case, the duality of Fractional Sobolev and Hardy-Littlewood-Sobolev inequalities are greatly related to a fast diffusion equation and use that diffusion flow to compare the remainder terms of the two dual inequalities for $n\ge 5$. Later Jin and Xiong \cite{jx13} and Dobeault and Jankowiak \cite{dj14} extended the result respectively to the case $s\in(0,2), n\ge 2, n>2s$ and the case $s=2, n\ge 3$. Actually, we have for $s=2, n\ge 3, p=q'=\frac{2n}{n+s}$,
there exists a positive constant $\alpha$, s.t., $~\forall~ 0\le  f\in \dot{H}^{s/2}$,
\beq\label{pr1-e} \norm{f}_{L^q}^{2(q-2)}\bra{\norm{f}_{\dot{H}^{s/2}}^2-C_{sharp} \norm{f}_{L^q}^2} \ge \alpha \bra{\norm{f^{q/p}}_{L^p}^2-C_{sharp}\norm{f^{q/p}}_{\dot{H}^{-s/2}}^2}.\eeq
Similar result for limit case --- Beckner-Onofri and Logarithmic Hardy-Littlewood-Sobolev inequalities are also given in \cite{dj14} in case $s=n=2$. We state on the sphere $\bb{S}^n$ and denote $\dashint =\frac{1}{|\bb{S}^n|}\int_{\bb{S}^n}$, then there exists a positive constant $\alpha$, s.t., if $\dashint e^f=1$, then
\beq\label{pr2-e} \frac{1}{2n!}\dashint f A_nf- \log \dashint e^{f-\dashint f} \ge \alpha \bra{\dashint e^f f-\frac{n!}{2}\dashint e^f A_n^{-1}e^f
},\eeq where $A_n$ is a spherical picture of $(-\Delta)^{n/2}$, 
meaning that $A_n f=\bra{|J_{\ca{S}}|^{-1}(-\Delta)^{n/2}(f\circ \ca{S})}\circ \ca{S}^{-1}$ with spectrum $j(j+1)\ldots(j+n-1)$ on the spherical harmonic subspace $\ca{H}_j$, which is injective when restricted on the image space with fundamental solution $-\frac{2}{(n-1)!}\log|\zeta-\eta|$ (note the kernel of $A_n$ is $\ca{H}_0$ and here $A_n^{-1}$ is interpreted after projection). The authors also gave some bounds about the sharp proportional constants and
we naturally think these sophisticated results about (\ref{pr1-e}) and (\ref{pr2-e}) can also be extended to \emph{all} fractions both on Eucliean spaces and Heisenberg groups. Indeed, later, after our doing for analogues on the Heisenberg group, we find corresponding results on $\bb{R}^n$ (about \ref{pr1-e} and \ref{pr2-e}) were also proved independently by Jankowiak and Hoang Nguyen in \cite{jn14} for all fractions $0<s<n, n\ge 2$.

On the Heisenberg group $\bb{H}^n$, the sharp fractional Sobolev-type (and its dual Hardy-Littlewood-Sobolev-type) inequalities were obtained recently in \cite{fl12} for all fractions. The extension to other Iwasawa-type groups were given by Christ and us in \cite{clz14,clz14-} for partial range of fractions. Some limit cases were also given in above papers and \cite{bfm13}. The main purpose of this note is to get, for all fraction $s$, some analogus results to Euclidean inequalities (\ref{sfs-e},\ref{pr1-e},\ref{pr2-e}) obtained in \cite{cfw13,dj14}, for conformally-invariant Sobolev-type inequalities on the Heisenberg group and more general Iwasawa-type groups. As usual, it's natural for us to carry out more conveniently the proof of the problems in the framework of spheres, where constant functions are extremizers for related inequalities. For the stability, we first give a local estimate and then use the recovery of compactness lemma to get a global one by contradiction. For the fractional Sobolev inequality, in the local neighborhood domain of $\ca{M}$ --- the submanifold of extremizers, with additional condition $d(f,\ca{M})<\norm{f}_{s/2}$, we can first use the conformal symmetry to assume the nearest point to be constant function 1 on the sphere and write the function to be $f=1+\varphi$, with $\varphi$ in the normal subspace, then the Taylor expansion tells us that we only need to compute and estimate the second variation around 1, which is observed to be positive definite on the normal subspace. For the Hardy-Littlewood-Sobolev inequality, we borrow a sophisticated local control of the functional aroud 1, due to \emph{Christ's lemma} for general $(p,q)$-functional with $p<2\le q$, to make up for the failure of Taylor expansion on $L^p$-distance.
To compare the two dual remainder terms, we first use the idea of \emph{completion of square method} to get a global proportional bound and then consider a local case by variational expansion, which in other word, give a upper bound of the best constant in the global proportional inequality.

\section{Fractional Sobolev and Hardy-Littlewood-Sobolev Inequalities}
Here, for simplicity, we consider on the Heisenberg group $\bb{H}^n$ the fractional Sobolev (Folland-Stein) inequality and its dual Hardy-Littlewood-Sobolev inquality, which both contain the intrinsic conformal invariance.

\emph{The Heisenberg Group.}
We identify the Heisenberg group $\bb{H}^n$ with its Lie algebra $\bb{C}^n\times\bb{R}$ endowed with group law
$uu'=(z,t)(z',t')=(z+z',t+t'+2\im z\cdot\overline{z'})$, for group elements $u=(z,t), u'=(z',t')$, where $z,z'\in \bb{C}^n, t,t'\in \bb{R}, z\cdot\overline{z'}=\sum_{j=1}^nz_j\overline{z_j'}$. The dilation for $\delta>0$ is $\delta: u=(z,t)\mapsto \delta(u)=(\delta z,\delta^2 t)$, and we denote the related homogeneous dimension by $Q=2n+2$, the homogeneous norm by $|u|=(|z|^4+|t|^2)^{1/4}$.
The left invariant vector fields, which coincide respectively with $\set{\frac{\partial}{\partial x_j}, \frac{\partial}{\partial y_j}, \frac{\partial}{\partial t}}_{j=1}^n$ at origin point, are given by
\[X_j=\frac{\partial}{\partial x_j}+2y_j\frac{\partial}{\partial t}, \quad Y_j=\frac{\partial}{\partial y_j}-2x_j\frac{\partial}{\partial t}, \quad T=\frac{\partial}{\partial t}.\] The \emph{sublaplacian} is a second order left invariant differential operator given by
\[\ca{L}=-\frac{1}{4}(X_j^2+Y_j^2),\] which is hypoelliptic from a celebrated H\"{o}rmander's theorem, and we recall that this essentially self-adjoint positive operator does not depend on the
choice of an orthonormal basis. An explicit computation gives the following formula
\[\ca{L}=-\frac{1}{4}\Delta_z-|z|^2\frac{\partial^2}{\partial t^2}+\frac{1}{2} N\frac{\partial}{\partial t},\] where \[\Delta_z=\sum_{j=1}^n\bra{\frac{\partial^2}{\partial x_j^2}+\frac{\partial^2}{\partial y_j^2}}, \quad N=\sum_{j=1}^n \bra{x_j\frac{\partial}{\partial y_j}-y_j\frac{\partial}{\partial x_j}}\] are respectively the standard Laplacian and  corresponding rotation operator on $\bb{C}^n$.
Using boundary Cayley transform on $\bb{H}^n$, a generalization of stereographic transform on Euclidean space $\bb{R}^n$, $\ca{C}: \bb{H}^n\rightarrow \bb{S}^{2n+1} \setminus \{o\}$ with $o$ being the south pole $(0,\ldots,0,-1)$, defined by
 \[u=(z,t)\mapsto \zeta=(\zeta',\zeta_{n+1})=\left(\frac{2z}{1+|z|^2-it},\frac{1-|z|^2+it}{1+|z|^2-it}\right),\] we can identify the Heisenberg group with the complex sphere. The Jacobian of the Cayley transform is
 \[|J_\ca{C}|=2^{Q-1}\bra{(1+|z|^2)^2+|t|^2}^{-Q/2}=2^{-1}|1+\zeta_{n+1}|^Q.\]
 A similar sublaplacian on $\s$ is defined from $\ca{L}$ by the Cayley transform and explicitly to be
 \[\ca{L}'=-\sum_{j=1}^{n+1} \frac{\partial^2}{\partial \zeta_j \partial \overline{\zeta_j}}+\sum_{j,k=1}^{n+1}\zeta_j\overline{\zeta_k}\frac{\partial^2}{\partial \zeta_j \partial \overline{\zeta_k}}+\frac{n}{2}\sum_{j=1}^{n+1}\bra{\zeta_j\frac{\partial}{\partial \zeta_j}+\overline{\zeta_j}\frac{\partial}{\partial \overline{\zeta_j}}},\]
 and the Geller-type \emph{conformal sublaplacian} is defined to be \[\ca{D}=\ca{L}'+\frac{n^2}{4}.\]
The ``sublaplacians" above play as important counterparts of positive Laplacian $-\Delta$ on Euclidean space $\bb{R}^n$ and there is an important relation between $\ca{L}$ and $\ca{D}$,
\[\ca{L}\bra{(2|J_\ca{C}|)^{\frac{Q-2}{2Q}}f\circ \ca{C}}=(2|J_{\ca{C}}|)^{\frac{Q+2}{2Q}}(\ca{D}f)\circ \ca{C}.\]
It's well known that the fundamental solutions of $\ca{L}$ and $\ca{D}$ are multiples of $(2-Q)$-power of distance functions:
\[\ca{L}^{-1}=\frac{2^{n-2}\Gamma^2(\frac{n}{2})}{\pi^{n+1}}\abs{u^{-1}v}^{2-Q}, ~ \ca{D}^{-1}=\frac{\Gamma^2(\frac{n}{2})}{2\pi^{n+1}}\abs{1-\zeta\cdot\bar{\eta}}^{\frac{2-Q}{2}}.\]

\emph{Fractional Sobolev Inequality.} The sharp fractional Sobolev (FS) inequality characterises the embedding $H^{s/2}\hookrightarrow L^q$, which generalizes a Jerison-Lee inequality (for $s=2$): \quad
 $\forall~ 0<s<Q, q=\frac{2Q}{Q-s}$ and $f\in H^{s/2}$, we have sharp inequality
\beq\label{fs}\norm{f}_{s/2}^2\ge C |f|_q^2,\eeq
where \beq\label{constant}C=
\bra{\frac{4\pi^{Q/2}}{n!}}^{s/Q}\frac{\Gamma^2(\frac{Q+s}{4})}{\Gamma^2(\frac{Q-s}{4})}\eeq is the best constant and will be fixed in this note. Here, for simplicity, we use $|\cdot|_q$ to denote the Lebesgue norm $\norm{\cdot}_{L^q}$, and $H^{s/2}$ is the fractional Sobolev space, which is the completion of $C_0^\infty$ w.r.t. the Sobolev norm $\norm{f}_{s/2}^2=\ip{f,f}_{s/2}=\ip{f,\ca{L}_sf}=\int f\ca{L}_sf$, with $\ca{L}_s$ being an \emph{intertwining operator} for complementary series representations of $SU(n+1,1)$. We have the following characterization of $\ca{L}_s$:
\[\ca{L}_s=|2T|^{s/2}\frac{\Gamma(\ca{L}|2T|^{-1}+\frac{2+s}{4})}{\Gamma(\ca{L}|2T|^{-1}+\frac{2-s}{4})}, \quad \ca{L}_s^{-1}=\frac{2^{n-1-s/2}\Gamma^2(\frac{Q-s}{4})}{\pi^{n+1}\Gamma(\frac{s}{2})}|u|^{s-Q}.\]
Note that $\norm{\cdot}_{s/2}$ is also equivalent to $|(I+\ca{L})^{s/2}|_2$, where $\ca{L}=\ca{L}_2$ is the sublaplacian --- the 2-order left invariant differential operator defined above. The FS inequality (\ref{fs}) is invariant under the conformal action $\sigma: f\mapsto \sigma(f)= f\circ\sigma |J_\sigma|^{1/q}$, where $\sigma$ is any conformal transformation on $\bb{H}^n$ with Jacobian determinant $|J_\sigma|$.
Note here $|\cdot|_q$ is obviously invariant while the inner product $\ip{\cdot,\cdot}_{s/2}$ is also invariant under the conformal actions, i.e., $\ip{\sigma(f),\sigma(g)}_{s/2}=\ip{f,g}_{s/2}$ considering the intertwining property of $\ca{L}_s$: \[|J_\sigma|^{\frac{Q+s}{2Q}}(\ca{L}_sf)\circ \sigma=\ca{L}_s\bra{|J_\sigma|^{\frac{Q-s}{2Q}}f\circ \sigma}
.\]
Using boundary Cayley transform $\ca{C}: \bb{H}^n\rightarrow \bb{S}^{2n+1} \setminus \{o\}$, 
we can move the inequality onto the complex sphere $\bb{S}^{2n+1}\subset \bb{C}^{n+1}$:
\beq\label{fs-}\norm{f}_*^2\ge C |f|_q^2,\eeq where $\norm{\cdot}_*$ is the norm induced by the inner product $\ip{f,g}_*=\ip{\ca{C}^{-1}(f), \ca{C}^{-1}(g)}_{s/2}$ and the two FS inequalities (\ref{fs})and (\ref{fs-}) are equivalent under the conformal correspondence of functions on the group and complex sphere $\ca{C}: f \mapsto \ca{C}(f)= f\circ \ca{C}|J_\ca{C}|^{1/q}$. This inner product $\ip{\cdot, \cdot}_*$ is a quadratic form involving an intertwining operator $\ca{A}_s$ on $\bb{S}^{2n+1}$, a spherical picture of $\ca{L}_s$. Acutally, $\ca{A}_s$ is uniquely given by the relation
\[\ca{L}_s\bra{|J_\ca{C}|^{\frac{Q-s}{2Q}}f\circ \ca{C}}=|J_{\ca{C}}|^{\frac{Q+s}{2Q}}(\ca{A}_sf)\circ \ca{C},\] and more precisely characterized by the bispherical harmonic decomposition or its fundamental solution as:
\beq\label{eig} \lambda_{j,k}=\ca{A}_s|_{\ca{H}_{j,k}}=2^{s/Q}\frac{\Gamma(j+\frac{Q+s}{4})\Gamma(k+\frac{Q+s}{4})}{\Gamma(j+\frac{Q-s}{4})\Gamma(k+\frac{Q-s}{4})},
\quad \ca{A}_s^{-1}=\frac{2^{-1-s/Q}\Gamma^2(\frac{Q-s}{4})}{\pi^{n+1}\Gamma(\frac{s}{2})}\abs{1-\zeta\cdot\bar{\eta}}^{\frac{2-Q}{2}},
\eeq
where $\ca{H}_{j,k}$ is the bispherical harmonic subspace spanned by harmonic polynomials of degree $j,k$ respectively on $\zeta$ and $\bar{\zeta}$ and we have the irreducible decomposition $L^2(\s)=\oplus_{j,k=1}^\infty \ca{H}_{j,k}$ (also $=\oplus_{j=1}^\infty \ca{H}_{j}$, where $\ca{H}_j$ is the classical real spherical harmonic subspace spanned by harmonic polynomials on real variables of degree $j$). There is also an intertwining relation for $\ca{A}_s$, just like that for $\ca{L}_s$, which, modula a constant, uniquely determines the operator. Note here $\ca{A}_2=2^{2/Q}\ca{D}$, where $\ca{D}$ is the Geller-type conformal sublaplacian and $\ca{D}|_{\ca{H}_{j,k}}=(j+\frac{n}{2})(k+\frac{n}{2})$.
We will use $H^*$ for the $\frac{s}{2}$-order Sobolev space on $\bb{S}^{2n+1}$, endowed with the norm $\norm{\cdot}_*$.
The conformal symmetry group of the complex spherical picture FS inequality (\ref{fs-}) consists of
\beq\label{inv}\tau: f\rightarrow \tau(f)=f\circ \tau |J_\tau|^{1/q}\eeq
 for any conformal transformation $\tau=\ca{C}\circ \sigma \circ \ca{C}^{-1}$ on $S^{2n+1}$, induced by the Cayley transform from any conformal transformation $\sigma$ on $\bb{H}^n$. We also have the invariance of the inner product $\ip{\cdot, \cdot}_*$ under the conformal actions.

It was proved in \cite{fl12} (in a dual reformulation) that all real-valued extremizers of the sharp FS inequality (\ref{fs}) are right all functions of the form $cH\big(\delta(u\cdot)\big)$, where $c\in\bb{R}\setminus \{0\}, \delta>0, u\in \bb{H}^n$, and \[H=\left((1+|z|^2)^2+|t|^2\right)^{-\frac{Q-s}{4}}=2^{-\frac{(Q-1)(Q-s)}{2Q}}|J_\ca{C}|^{1/q},\] i.e., (\ref{fs}) reaches equality only on $cH$, up to the conformal invariant actions. Similarly, all real-valued extremizers of (\ref{fs-}) are right given by $c\abs{1-\xi\cdot\bar{\zeta}}^{-\frac{Q-s}{2}}$, where $c\in \bb{R}\setminus\{0\}, \xi \in \bb{C}^{n+1}$ and $|\xi|<1$,  i.e., (\ref{fs-}) reaches equality only on constant functions, up to the conformal symmetry group. Take $\ca{M}, \ca{M}_*$ respectively to be the $(Q+1)$-dimensional smooth submanifolds of all real-valued extremizers for (\ref{fs}) and (\ref{fs-})and define naturally the distances
$d(f,\ca{M}):=\inf\{\norm{f-g}_{s/2}: g\in \ca{M}\}, ~d(f,\ca{M}_*):=\inf\{\norm{f-g}_*: g\in \ca{M}_*\}$. Then we have the following theorem characterising the stability of extremizers.
\bthm\label{t-fs} \emph{(Stability for FS)}
Let $0<s<Q=2n+2, q=\frac{2Q}{Q-s}$, then there exists $\alpha>0$ only depending on the dimension $Q$ and fraction $s$, s.t.,
\begin{align}
d^2(f,\ca{M})\ge \norm{f}_{s/2}^2-C|f|_q^2\ge \alpha~ d^2(f,\ca{M}),&\quad \forall f\in H^{s/2}, \label{f-fs}\\
d^2(f,\ca{M}_*)\ge \norm{f}_*^2-C|f|_q^2\ge \alpha~ d^2(f,\ca{M}_*),&\quad \forall f\in H^*,\label{f-fs-}
\end{align} and if $d(f,\m)$ or $d(f,\ca{M}_*)>0$, the left inequality is strict.
\ethm

\emph{Hardy-Littlewood-Sobolev Inequality.} A duality argument writes the FS inequality (\ref{fs}) into the (diagonal) Hardy-Littlewood-Sobolev (HLS) inequality, which gives the $(p,p')$-boundness of the fractional integral operator $|u|^{-\lambda}*$: \quad $\forall~ 0<\lambda<Q$, ~$p=\frac{2Q}{2Q-\lambda}$ and $f,g\in L^p$, we have
\[ \abs{|u|^{-\lambda}*f}_{p'}\lesssim |f|_p,  \quad \text{or equivalently} \quad \abs{\iint_{\bb{H}^n\times\bb{H}^n} \frac{f(u)g(v)}{|u^{-1}v|^\lambda}dudv} \lesssim |f|_p|g|_p,\] and by the boundary Cayley transform, there is an equivalent HLS inequality on the complex sphere $\s$,
\[\abs{\iint_{\s\times\s} \frac{f(\zeta)g(\eta)}{\abs{1-\zeta\cdot\bar{\eta}}^{\lambda/2}}d\zeta d\eta} \lesssim |f|_p|g|_p.\]
Modulo a constant multiple (see from the fundamental solution formulas of $\ca{L}_s$ and $\ca{A}_s$), we write above HLS inequalities in the following sharp form: \quad $\forall~ 0<s<Q, ~p=\frac{2Q}{Q+s}$, ~$f\in L^p$ (on $\bb{H}^n$ or $\s$ depending on the corresponding context),
\beq\label{hls}\norm{f}_{-s/2}^2 \le C^{-1} |f|_p^2, \quad \norm{f}_{-*}^2 \le C^{-1} |f|_p^2,\eeq where $\norm{\cdot}_{-s/2}$ is the negative fractional Sobolev norm on $\h$ under the meaning of $\norm{\cdot}_{-s/2}^2=\ip{\cdot,\ca{L}_s^{-1}\cdot}$ and $\norm{\cdot}_{-*}$ is the correspondence on $\s$ induced from the inner product $\ip{\cdot,\cdot}_{-*}=\ip{\cdot,\ca{A}_s^{-1}\cdot}$. The two HLS inequalities (\ref{hls}) are respectively invariant under the conformal action $\sigma: f\mapsto f\circ \sigma |J_\sigma|^{1/p}$ and $\tau: f\mapsto f\circ \tau |J_\tau|^{1/p}$, where $\sigma$ and $\tau$ are respectively the conformal transformation on the group and sphere. This can be seen from the conformal invariance of the FS inequalities (\ref{fs},\ref{fs-}) by dual argument or derived directly from the following formulas:
\[\abs{\sigma(u)^{-1}\sigma(v)}=\abs{u^{-1}v}\abs{J_\sigma(u) J_\sigma(v)}^{\frac{1}{2Q}}, \abs{1-\tau(\zeta)\cdot \overline{\tau(\eta)}}^{1/2}=\abs{1-\zeta\cdot\bar{\eta}}^{1/2}\abs{J_\tau(\zeta) J_\tau(\eta)}^{\frac{1}{2Q}}.\]
Take the submanifold of real-valued extremizers for HLS inequalites on the group and sphere respectively to be
\begin{align*}
\ca{M}_-=\set{c\bra{\abs{J_\ca{C}}^{1/p}}\big(\delta(u\cdot)\big): c\in \bb{R}\setminus\{0\}, \delta>0, u\in \bb{H}^n}, \qquad ~ \qquad \ ~ \\
\ca{M}_{-*}=\Big\{c|J_\tau|^{1/p}=c\abs{1-\xi\cdot\bar{\zeta}}^{-\frac{Q+s}{4}}: c\in \bb{R}\setminus\{0\}, \xi\in\bb{C}^{n+1}, |\xi|<1, ~ \  \\
\tau ~\text{is a conformal transformation on}~\s \Big\}.
\end{align*}
Define the $L^p$ distances to the extremizering submanifolds  to be $d_p(f,\m_-)=\inf\{|f-g|_p: g\in \m_-\}$, and $d_p(f,\m_{-*})=\inf\{|f-g|_p: g\in \m_{-*}\}$ and denote $H^{-s/2}, H^{-*}$ for the associated negative fractional Sobolev spaces on the group and sphere. We then have a similar stability result for the HLS inequalities (\ref{hls}).
\bthm\label{t-hls} \emph{(Stability for HLS)}
Let $0<s<Q=2n+2, ~p=\frac{2Q}{Q+s}$, then
there exist $\alpha_0,\alpha_1>0$ only depending on the dimension $Q$ and fraction $s$, s.t., $\forall~ 0\not\equiv f\in L^p$ (on the group or sphere depending on the context)
\begin{align}
 \alpha_1~ \frac{d_p(f,\ca{M}_-)}{|f|_p} \ge C^{-1}-\frac{\norm{f}_{-s/2}^2}{|f|_p^2} \ge \alpha_0~ \frac{d_p^2(f,\ca{M}_-)}{|f|_p^2}, \label{f-hls}\\
\alpha_1~ \frac{d_p(f,\ca{M}_{-*})}{|f|_p} \ge C^{-1}-\frac{\norm{f}_{-*}^2}{|f|_p^2} \ge \alpha_0~ \frac{d_p^2(f,\ca{M}_{-*})}{|f|_p^2}.\label{f-hls-}
\end{align}
\ethm
We remark a difference from the problem for the FS inequality. For the upper bound, we can only do for first power as we don't have any orthogonality, and still don't know about the square distance. This global bound, a stability, was still got from some local estimate. For the lower local bound, the Taylor expansion fails as the exponent $p$ is strictly less than $2$ and however, we can borrow one useful lemma of Christ's, who first notice the failure of $C^2(L^p)$ of the functional as $p<2$ and use very sophisticated estimates to get a bound for the functional using second variation of suitably chosen small truncation. See details in subsection \ref{proof of shls}. 

\emph{Relation Between Remainder Terms of FS and HLS Inequalities.} Now we want to compare the dual remainder terms of FS and HLS inequalities by a constant multiple, both globally and locally. The similar problem on Euclidean spaces has been studied in \cite{jx13,dj14} and is now extended to the Heisenberg group for total range of $s$.
\bthm\label{t-relation} \emph{(Dual remainder terms inequality)}
About remainder terms of the sharp FS and HLS inequalities on $\s$, we have the following estimates (an equivalence on $\h$ exists):
\benu[$1)$]
\item \emph{(Global estimate)}
\beq\label{global bound}|f|_q^{2(q-2)}(\norm{f}_*^2-C|f|_q^2) \ge  C\bra{\abs{f^{q/p}}_p^2-C\norm{f^{q/p}}_{-*}^2}, \qquad \forall~0\le f\in H^*.\eeq
\item\emph{(Local estimate)}
\beq\label{local bound}\liminf_{0<d(f,\m_*)\rightarrow 0, ~0\le f\in H^*,  ~d(f,\m_*)<\norm{f}_*}\frac{|f|_q^{2(q-2)}(\norm{f}_*^2-C|f|_q^2)}{\abs{f^{q/p}}_p^2-C\norm{f^{q/p}}_{-*}^2} = \frac{Q+4+s}{Q+4-s} ~C.
\eeq
\eenu
\ethm

The above results can be extended to all Iwasawa-type groups. We list on the sphere the results without proof as they are absolutely similar to the Heisenberg group case. First, we introduce something about groups of Iwasawa-type, which can be seen as the nilpotent part of the Iwasawa decomposition of a semisimple Lie group of rank one. It's one kind of groups of Heisenberg-type, satisfying additional $J_2$-condition. A two-step nilpotent Lie algebra $\fr{g}$ with center $\fr{z}$, endowed with an inner product $\ip{\cdot,\cdot}$, for which $\fr{g}=\fr{z}^\perp\oplus\fr{z}$, is said to be of Heisenberg-type, if the Lie structure satisfies the following condition: for any element $|t|=1$ in the center $\fr{z}$, an endomorphism $J_t$ defined on $\fr{z}^\perp$ by $\ip{J_tz,z'}=\ip{t,[z,z']}$ for any $z,z'\in \fr{z}^\perp$, is an orthogonal map. A simply connected connected Lie group is said to be of  Heisenberg type if its Lie algebra is of Heisenberg type.  A $J_2$-condition is defined to be: for any $t,t'\in \fr{z}$ satisfying $\ip{t,t'}=0$, there exists a $t''\in \fr{z}$, s.t., $J_{t''}=J_tJ_t'$. It was well known that any group of Heisenberg-type is of Iwasawa-type if and only if the Lie structure satisfies the $J_2$-condition. Besides, from a geometrical point of view, the Iwasawa groups can be seen as the nilpotent part in the Iwasawa decomposition of
the isometry group of the associated non-compact Riemannian symmetric spaces of rank one and finally Iwasawa-type groups are identified to four cases having center of dimension 0,1,3,7. We write isometrically the group in a unifying form to be $\bb{K}^n\times \im \bb{K}$ for appropriate $n$, where $\bb{K}$ is one of the four real division algebra $\bb{R}, \bb{C},\bb{H}, \bb{O}$. For above results and FS, HLS and related inequalities in the later two cases, see \cite{clz14,clz14-} and references there. If we still use the same notation as those on the Heisenberg group and denote $m$ to be the dimension of center, we can state the following result.
\bthm \emph{(For general Iwasawa-type group)}
We have analogous sharp FS and HLS inequalities on the Iwasawa-type group for fractions of partial range $s < Q-4[\frac{m}{2}], \lambda > 4[\frac{m}{2}]$, where $m$ is the dimension of center and $Q$ is the homogeneous dimension. Moreover, Theorem \ref{t-fs},\ref{t-hls} and \ref{t-relation} still holds similarly on the Iwasawa-type group for the corresponding inequality with exponent of partial range. More precisely, we have the following properties for the remainder terms: let $p=q'=\frac{2Q}{Q+s}$,
\benu[$1)$]
\item \emph{(Stability for FS and HLS)}  There exist two positive constants $\alpha_0, \alpha_1$ only depending on $Q$ and $s$, s.t.,
\begin{align*}
d^2(f,\ca{M}_*)\ge \norm{f}_*^2-C|f|_q^2 \ge \alpha_0~ d^2(f,\ca{M}_*), \quad \forall f\in H^*,\\
\alpha_1~ |f|_p d_p(f,\ca{M}_{-*}) \ge |f|_p^2 - C\norm{f}_{-*}^2 \ge \alpha_0~ d_p^2(f,\ca{M}_{-*}), \quad \forall f\in L^p,
\end{align*} where $C,\m,\m_*$ is respectively the sharp constant and extremizering submanifolds, computed in \emph{\cite{lie83,fl12,clz14,clz14-}}, and $d(f,\ca{M}_*)$, $d_p(f,\m_{-*})$ is respectively the Sobolev and Lebesgue distance to the submanifold.
\item\emph{(Dual remainder terms inequality)}
About remainder terms of the sharp FS and HLS inequalities on the sphere, we have 
\benu[$a)$]
\item\emph{(Global estimate)}
\[|f|_q^{2(q-2)}(\norm{f}_*^2-C|f|_q^2) \ge  C\bra{\abs{f^{q/p}}_p^2-C\norm{f^{q/p}}_{-*}^2}, \qquad \forall~ 0\le f\in H^*.\]
\item\emph{(Local estimate)}
\[\liminf_{0<d(f,\m_*)\rightarrow 0,~0\le f\in H^*, ~d(f,\m_*)<\norm{f}_*}\frac{|f|_q^{2(q-2)}(\norm{f}_*^2-C|f|_q^2)}{\abs{f^{q/p}}_p^2-C\norm{f^{q/p}}_{-*}^2} = \frac{Q+2+2\,\emph{sign}m+s}{Q+2+2\,\emph{sign}m-s} ~C.
\]
\eenu
\eenu
\ethm

We add a note here about the sharp proportional constant between the two dual remainder terms, which is still open even for the Eucliean case. The global estimate tells that the sharp constant, denoted by $C_{sharp}$, is bigger than $C$, i.e., $C_{sharp}\ge C$, while the local estimate tells $C_{sharp}\le \frac{Q+4+s}{Q+4-s} ~C$. However, we expect a fast diffusion method like in \cite{dj14} may improve a bit the lower bound strictly to be $C_{sharp}>C$.
\section{Limit Case --- Beckner-Onofri and Logarithmic Hardy-Littlewood-Sobolev Inequalities.}
Now, as in \cite{fl12} and \cite{bfm13}, we can consider on the Heisenberg group the endpoint limit case ($s=Q, \lambda=0$), using the functional differential argument. We use $\dashint$ to denote the mean integral on the sphere, i.e., $\dashint f= |\s|^{-1}\int_{\s} f$.

\emph{Beckner-Onofri Inequality.}
For the limit case $s=Q$, the generalization of sharp Beckner-Onofri (BO) inequality on the Heisenberg group and complex sphere was first given by \cite{bfm13}. For simplicity, we only deal on the sphere, which involves the \emph{conditional intertwining operator} $\ca{A}_{Q}'$, a kind of differentiation of $\ca{A}_s (0<s<Q)$ at $Q$. The operator, defined on $H^*\cap\ca{P}$, where $H^*$, in abuse of notation, is the $\frac{Q}{2}$-order Sobolev space and $\ca{P}=\oplus_{jk=0} \ca{H}_{j,k}$ ($\bb{R}\ca{P}$) is the space of $L^2$ (real-valued) CR-pluriharmonic functions, is given by:
\beq\label{eig'} \lambda_j=\ca{A}_Q'|_{\ca{H}_{j,0}}=\ca{A}_Q'|_{\ca{H}_{0,j}}=\frac{\Gamma(j+n+1)}{\Gamma(j)}=j(j+1)\ldots(j+n).\eeq
If in this subsection, we denote (a constant difference from (\ref{eig}) in last section) \beq\label{mod}\lambda_{j,k}=\ca{A}_s|_{\ca{H}_{j,k}}=\frac{\Gamma(j+\frac{Q+s}{4})\Gamma(k+\frac{Q+s}{4})}{\Gamma(j+\frac{Q-s}{4})\Gamma(k+\frac{Q-s}{4})},\eeq then
\beq\label{dif}\ca{A}_Q'=-\frac{4}{n!} \frac{d}{ds}\Big|_{s=Q} \ca{A}_s= \lim_{s\rightarrow Q}\frac{4}{n!}\frac{1}{Q-s}\ca{A}_s.\eeq
The sharp BO inequality states that: \quad $\forall f\in H^*\cap\bb{R}\ca{P}$ (with zero mean), we have
\beq\label{bo}\frac{1}{2(n+1)!}\dashint f\ca{A}_Q'f  \ge \log\dashint e^{(f-\dashint f)}\eeq (without the mean term on the right side). Obviously, the BO inequality (\ref{bo}) is invariant up to constant translations, and besides, there is also a big conformal invariance action group for this inequality, \beq\label{inv-}\tau: f\mapsto \tau(f)= f\circ \tau +\log |J_\tau|.\eeq We can easily see the right side exponential integral $\dashint e^f$ is invariant under the conformal action, while
we can use the limit differentiation argument and conformal invariance property of FS inequality to get the whole invariance for the BO inequality (\ref{bo}): from (\ref{constant}), we define the FS functional of (\ref{fs-}) modified by (\ref{mod}) to be
\[I(g)=\dashint g \ca{A}_sg-\frac{\Gamma^2(\frac{Q+s}{4})}{\Gamma^2(\frac{Q-s}{4})}\bra{\dashint|g|^q}^{2/q},\] then it satisfies the following formulas: take $\lambda=Q-s\rightarrow 0,~ \forall~ f\in H^*\cap \bb{R}\ca{P}$,\\
(1) If $\dashint f=0$, we have
\begin{align*}
I(1+\lambda f)=&~\ip{1,\ca{A}_s1}+\lambda^2\ip{f,\ca{A}_sf}-\bra{\bra{\frac{\lambda}{4n!}}^2+o(1)}\bra{\dashint e^{2Qf}+o(1)}^{\lambda/Q}\\
=&~ \frac{n!}{4}\ip{f,\ca{A}_Q'f} \lambda^3-\bra{\bra{\frac{n!}{4}}^2\lambda^2+o(1)}\bra{\bra{\dashint e^{2Qf}+o(1)}^{\lambda/Q}-1}\\
=&~ \frac{n!}{4}\ip{f,\ca{A}_Q'f} \lambda^3-\bra{\bra{\frac{n!}{4}}^2\lambda^2+o(1)}\bra{\frac{1}{Q}\log\dashint e^{2Qf}\lambda+o(\lambda)}\\
=&~ \frac{1}{Q}\bra{\frac{n!}{4}}^2 \bra{\frac{4Q}{n!}\dashint f\ca{A}_Q'f-\log \dashint e^{2Qf}}\lambda^3+o(\lambda^3),
\end{align*} where we use notation $\ip{f,g}=\dashint f g$.\\
(2) for general $f$, we have
\[ I\bra{1+\frac{\lambda}{2Q}f}= \frac{1}{Q}\bra{\frac{n!}{4}}^2 \bra{\frac{1}{2(n+1)!}\dashint f\ca{A}_Q'f-\log \dashint e^{f-\dashint f}}\lambda^3+o(\lambda^3).\]
(3) from the invariance of the FS inequality functional under the corresponding conformal action (\ref{inv}), we have
\begin{align*}
I\bra{1+\frac{\lambda}{2Q}f}=&~I\bra{\bra{1+\frac{\lambda}{2Q}f}\circ\tau ~|J_\tau|^{1/q}}\\
=&~I\bra{1+\frac{\lambda}{2Q}(f\circ\tau+\log |J_\tau|)+o(\lambda)}\\
=&~I\bra{1+\frac{\lambda}{2Q}(f\circ\tau+\log |J_\tau|)}+ o(\lambda^3).
\end{align*}
Then we see the invariance of (\ref{bo}) under (\ref{inv-}) from
\[\frac{1}{2(n+1)!}\dashint f\ca{A}_Q'f  - \log\dashint e^{f-\dashint f}=Q\bra{\frac{4}{n!}}^2 \frac{I\bra{1+\frac{\lambda}{2Q}f}}{\lambda^{3}}+o(1).\]
We know extremizers for BO inequality (\ref{bo}) are right all the functions of the following form $\log (c|J_\tau|)=\log c+Q\log \abs{1-\xi\cdot\bar{\zeta}}^{-1}$, where $c>0, \xi\in \bb{C}^{n+1}, |\xi|<1$ and $\tau$ is any conformal transformation on $\s$ ($c=1$ when adding zero mean).

For stability, it seems that we can't find a good distance norm, which is invariant under the conformal action, to characterize the relation between remainder terms and distance to $(Q+1)$-dimensional smooth submanifold of all real-valued extremizers. (The norm $\dashint f\ca{A}_Q'f$, and also the subspace of zero-mean functions in $\ca{P}$, are not invariant under the conformal action, while the functional $\frac{1}{2(n+1)!}\dashint f\ca{A}_Q'f-\dashint f$ is unable to characterize the distance).

\emph{Logarithmic Hardy-Littlewood-Sobolev Inequality.}
The generalized sharp Logarithmic Hardy-Littlewood-Sobolev (Log-HLS) Inequality on the Heisenberg group and complex sphere was first obtained by \cite{bfm13} and reproved in \cite{fl12}, which states that: \quad
for any nonnegative normalized $f\in L\log L, \dashint f=1$, we have
\[(n+1)\dashint\dashint \log \frac{1}{|1-\zeta\cdot\bar{\eta}|}f(\zeta)f(\eta)d\zeta d\eta\le \dashint f\log f.\] Because the fundamental solution of operator $\ca{A}_Q'$ restricted on the image subspace is given by \[\ca{A}_Q'^{-1}\big|_{f\in\ca{P},\dashint f=0}=\frac{1}{\pi^{n+1}}\log \frac{1}{|1-\zeta\cdot \bar{\eta}|},\] above Log-HLS inequality is equivalent to the following dual reformulation of the BO inequality (\ref{bo}):
\beq\label{loghls}\frac{(n+1)!}{2}\dashint(f-1)\ca{A}_Q'^{-1}\bra{P(f-1)}\le \dashint f\log f,\eeq where $P$ denotes the projection operator from $L^2$ onto $\ca{P}$ and for simplicity we often leave $P$ out and use directly $\ca{A}_Q'^{-1}$ for $\ca{A}_Q'^{-1}\circ P$. The Log-HLS inequality is invariant under the conformal action $f \mapsto f\circ \tau |J_\tau|$.

\emph{Relations Between Remainder Terms of BO and Log-HLS Inequality.}
Just as estimate for the FS and HLS remainder terms, we can get both a global and a local bound for the BO and Log-HLS remainder terms. Differentiation is used to get the global one while the local one is got directly from Taylor expansion and the property of variations.
\bthm \label{t-relation-} \emph{(Dual inequality for BO and Log-HLS)}
About remainder terms of the sharp BO and Log-HLS inequalities on $\s$, we have the following estimates (an equivalence on $\h$ exists):
\benu[$1)$]
\item\emph{(Global estimate)}
\begin{align}
\frac{1}{2(n+1)!}\dashint f\ca{A}_Q'f-\log \dashint e^{f-\dashint f} \ge
\frac{1}{\dashint e^f}\dashint e^ff - \frac{(n+1)!}{2(\dashint e^f)^2} \dashint\! \big(e^f-\dashint e^f\big)\ca{A}_Q'^{-1}\big(e^f-\dashint e^f\big)\nonumber\\
 -\log \dashint e^f, ~\qquad \forall f\in H^*\cap \bb{R}\ca{P}.\label{global-limit}
\end{align}
Note the invariance of left and right side up to constant translation $f\mapsto f-c$, then for any $f$, we can choose $c=\log \dashint e^f$, s.t., after constant translation, new $f$ satisfies $\dashint e^f=1$ and the inequality becomes
\[\frac{1}{2(n+1)!}\dashint f\ca{A}_Q'f-\log \dashint e^{f-\dashint f} \ge \dashint e^f f-\frac{(n+1)!}{2}\dashint\bra{e^f-1}\ca{A}_Q'^{-1}\bra{e^f-1}.\]
\item\emph{(Local estimate)}
\beq\label{local-limit}\liminf_{0<\dashint f\ca{A}_Q'f \rightarrow 0}^{f\in H^*\cap\ca{P}, f\perp \ca{H}_1} \frac{\frac{1}{2(n+1)!}\dashint f\ca{A}_Q'f-\log \dashint e^{f-\dashint f}}{\frac{1}{\dashint e^f}\dashint e^ff - \frac{(n+1)!}{2(\dashint e^f)^2}
\dashint \big(e^f-\dashint e^f\big)\ca{A}_Q'^{-1}\big(e^f-\dashint e^f\big) -\log \dashint e^f} \ge n+2.\eeq
\eenu
\ethm

On the two other Iwasawa-type groups, there is difficulty in getting the sharp FS and HLS inequality for big $s$ and small $\lambda$, as the eigenvalue can be negative. Therefore, endpoint case --- some putative sharp BO and Log-HLS inequalities are also still unknown.

\section{Proof of Main Results}
We only prove the second formula in Theorem \ref{t-fs} and \ref{t-hls}, recalling the equivalence of that two formulas on the group and sphere.
\subsection{Proof of FS Stability}
For Theorem \ref{t-fs}, we first need a local result.
\bpro\label{local} \emph{(Local stability for FS)}
We have the following local reminder term estimate: for all $f\in H^*$, satisfying $d(f,\ca{M}_*)<\norm{f}_*$,
\beq\label{f-local}d^2(f,\ca{M}_*)\ge \norm{f}_*^2-C|f|_q^2\ge \frac{2s}{Q+4+s} d^2(f,\ca{M}_*)+o\bra{d^2(f,\ca{M}_*)},\eeq and if $d(f,\ca{M}_*)>0$, the left inequality is strict.
\epro
\bpf
Consider any $f\in H^*$, from definition, there exists a function $g\in \overline{\ca{M}_*}$, s.t., $\norm{f-g}_*=d(f,\ca{M}_*)$. From condition $d(f,\ca{M}_*)<\norm{f}_*$, we know $g \not\equiv 0$, and there exist $c \in \bb{R}\setminus \{0\}, \xi\in \bb{C}^{n+1}, |\xi|<1$, and a conformal transformation $\tau$, s.t. $g=c\tau(1)=c|J_\tau|^{1/q}=c\abs{1-\xi\cdot \bar{\zeta}}^{-\frac{Q-s}{2}}$. So, from the invariance of $\ca{M}_*$ and above inequality (\ref{f-local}) under the conformal action, we may assume $f=1+\varphi$ with $\varphi \perp T_1 \ca{M}_*$ (normal to the tangent space at point --- function 1) under the inner product $\ip{\cdot, \cdot}_*$.  We take the difference functional $I(f)=\norm{f}_*^2-C|f|_q^2$, which is in $C^2(H^*)$, then for local lower bound, by Taylor expansion and the critical property of extremizer 1, we have
\beq\label{taylor}I(f)=\frac{1}{2}I''(1,\varphi)+o(\norm{\varphi}_*^2),\eeq
where
\begin{align*}
I'(1,\varphi)=\frac{d}{dt}\Big|_{t=0}I(1+t\varphi)= 2\ip{1,\varphi}_*-2C|1|_q^{2-q}\int \varphi \equiv 0.\\
I''(1,\varphi)=\frac{d^2}{dt^2}\Big|_{t=0} I(1+t\varphi)=2\norm{\varphi}_*^2-2(q-1)C|1|_q^{2-q}|\varphi|_2^2.
\end{align*}
More precisely, for the Taylor expansion (\ref{taylor}), we only need to prove \beq\label{taylor-}|1+\varphi|_q^2=|1|_q^2+(q-1)|1|_q^{2-q}|\varphi|_2^2+o(|\varphi|_q^2), \quad~\forall~ q>2, \varphi\in L^q, \int \varphi=0.\eeq
We note that the quantity $|1+\varphi|^q$ can't be expanded by Taylor series around 1 unless $|\varphi|$ is very small relative to 1, this gives the difficulty for the Taylor expansion of the difference (or $L^q$) functional. To prove (\ref{taylor-}), we use a truncation considering for $|\varphi|\le\delta$ and $|\varphi|>\delta$ respectively, setting $\delta=|\varphi|_q^{1-\epsilon}$ and $|\varphi|_q\rightarrow 0$,
\begin{align*}
\forall~|\varphi|\le \delta, \quad & |1+\varphi|^q=1+q\varphi+\frac{q(q-1)}{2}\varphi^2+\frac{q(q-1)(q-2)}{3!}|1+\theta\varphi|^{q-3}\varphi^3,\\
\forall~|\varphi|> \delta, \quad & \frac{\abs{|1+\varphi|^q-(1+q\varphi)}}{|\varphi|^q}\le  \bra{\frac{q(q-1)}{2}+\epsilon} \delta^{2-q}, ~\quad |\varphi|^2<\delta^{2-q}|\varphi|^q.
\end{align*}
then we have a uniform estimate for $q>2$,
\[\abs{|1+\varphi|^q-(1+q\varphi+\frac{q(q-1)}{2}\varphi^2)}\lesssim \delta\varphi^2+\delta^{2-q}|\varphi|^q,\]  and by integration, we have
\[\abs{|1+\varphi|_q^q-(|1|_q^q+\frac{q(q-1)}{2}|\varphi|_2^2)}\lesssim \delta|1|_q^{q-2}|\varphi|_q^2+\delta^{2-q}|\varphi|_q^q=o(|\varphi|_q^2),\] which gives the assertion (\ref{taylor-}) and therefore (\ref{taylor}) after raising to the power $2/q$ as
\[|1+\varphi|_q^2=\bra{|1|_q^q+\frac{q(q-1)}{2}|\varphi|_2^2+o(|\varphi|_q^2)}^{\frac{2}{q}}=|1|_q^2+\frac{2}{q}\frac{q(q-1)}{2}|1|_q^{2-q}|\varphi|_2^2+o(|\varphi|_q^2).\]
Now we return to the Taylor expansion of the difference functional and to prove the local lower bound.
Note that $T_1\ca{M}_*=Span\{1, \re\zeta_j, \im\zeta_j\}_{j=1}^{2n+2}=\ca{H}_0\oplus\ca{H}_1=\ca{H}_{0,0}\oplus\ca{H}_{1,0}\oplus\ca{H}_{0,1}$ and from bispherical harmonics expansion $\varphi=\sum_{j,k}\varphi_{j,k}$, where $\varphi_{j,k}\in \ca{H}_{j,k}$, we see
\[\frac{|\varphi|_2^2}{\norm{\varphi}_*^2}= \frac{\sum_{j+k\ge2}|\varphi_{j,k}|_2^2}{\sum_{j+k\ge2}\lambda_{j,k}|\varphi_{j,k}|_2^2},\] and
\begin{align*}
I(f)=&~\norm{\varphi}_*^2-(q-1)C\abs{\bb{S}^{2n+1}}^{\frac{2-q}{q}}|\varphi|_2^2+o(\norm{\varphi}_*^2)\\
=&~\norm{\varphi}_*^2\bra{1-(q-1)C\abs{\bb{S}^{2n+1}}^{\frac{2-q}{q}}\frac{|\varphi|_2^2}{\norm{\varphi}_*^2}+o(1)}.
\end{align*}
From the increasing of $\lambda_{j,k}$ on $j,k$ (see (\ref{eig}) and note $\lambda_{1,1}\ge \lambda_{2,0}$), we have $\frac{|\varphi|_2^2}{\norm{\varphi}_*^2} \le \lambda_{2,0}^{-1}$ and recalling the fact $\frac{\lambda_{1,0}}{C}=(q-1)\abs{\s}^{\frac{2-q}{q}}$ (see (\ref{constant}) and (\ref{eig})), we finally get
\begin{align*}
I(f)\ge&~ \norm{\varphi}_*^2 \bra{1-(q-1)C\abs{\bb{S}^{2n+1}}^{\frac{2-q}{q}}\lambda_{2,0}^{-1}+o(1)}\\
=&~ \norm{\varphi}_*^2\bra{\frac{2s}{Q+4+s}+o(1)}.
\end{align*} So, for all $s\in (0,Q)$, we give a positive local lower bound.
For the upper bound, we use $\ip{\varphi,1}_*=0$ and get
\[I(f)=\norm{\varphi}_*^2+\norm{1}_*^2-C|1+\varphi|_q^2 \le \norm{\varphi}_*^2,\] recalling $|1+\varphi|_q^2\ge \abs{\bb{S}^{2n+1}}^{2/q-1}|1+\varphi|_2^2\ge \abs{\bb{S}^{2n+1}}^{2/q}=|1|_q^2$, which reaches equality if and only if $\varphi\equiv 0$, or equivalently  $d(f,\m_*)=0$. The inequality (\ref{f-local}) and therefore Proposition \ref{local} is proved.
\epf

For global estimate of the stability, we also need the following lemma asserting the possibility of recovering compactness from the conformal symmetries. The lemma is proved in a technical $TT^*$ way for dual HLS inequality in \cite{fl12} and can also be derived from suitably adapted \emph{concentration compactness} argument, originated from Lions's work for HLS inequality on $\bb{R}^n$ \cite{lio85}. It's also natural to apply directly the adapted profile decomposition from G\'{e}rard's work \cite{ger98} on Euclidean spaces.
\blem\label{compact}\emph{(Recovery of compactness for FS)}
Let $(f_j)$ be an (non-vanishing) extremizing sequence of inequality \emph{(\ref{fs-})} (or its functional), i.e., $\frac{\norm{f_j}_*^2}{|f_j|_q^2}\xrightarrow{j\rightarrow\infty} C$, then  $\frac{d(f_j,\ca{M}_*)}{\norm{f_j}_*} \xrightarrow{j\rightarrow \infty} 0$.
\elem

We can now use the local stability --- Proposition \ref{local} and above recovery of compactness property --- lemma \ref{compact} to prove the global stability for FS inequality.\\

\emph{Proof of Theorem \ref{t-fs}.}
Core is the lower bound for (\ref{f-fs-}). For upper bound, as $0\in \overline{\ca{M}_*}$, we have $d(f,\ca{M}_*)\le \norm{f}_*$, then it's trivial from upper bound of (\ref{f-local}) in Proposition \ref{local}. For lower bound of (\ref{f-fs-}), we argue by contradiction (this works only for existence of a positive constant). If the constant $\alpha>0$ doesn't exist, then we can get a sequence $(f_j)$ satisfying $\frac{\norm{f_j}_*^2-C|f_j|_q^2}{d^2(f_j,\ca{M}_*)}\xrightarrow{j\rightarrow \infty} 0$ with $d(f_j,\ca{M}_*)>0$. From $\frac{d^2(f,\ca{M}_*)}{\norm{f}_*^2} \le 1$, we have $\frac{\norm{f_j}_*^2}{|f_j|_q^2} \xrightarrow{j\rightarrow\infty} C$, and then from Lemma \ref{compact}, we have $\frac{d(f_j,\ca{M}_*)}{\norm{f_j}_*} \xrightarrow{j\rightarrow \infty} 0$, which tells that $(f_j)$ is not only an extremizering sequence but also a local one to $\m_*$, whence implies from the local lower bound of (\ref{f-local}) in Proposition \ref{local} that, $\liminf_{j\rightarrow \infty}\frac{\norm{f_j}_*^2-C|f_j|_q^2}{d^2(f_j,\ca{M}_*)} \ge \frac{2s}{Q+4+s}>0.$ So we have got a contradiction and therefore proven the global theorem totally. \qed \\

\subsection{Proof of HLS Stability}\label{proof of shls}

We now still use the contradiction idea to prove Theorem \ref{t-hls}, the global stability of the HLS inequality. For (\ref{f-hls-}), first we need to consider locally, and recalling the failure of Taylor expansion of associated functional (see the second variation in \emph{a priori} Taylor expansion $|\varphi|_2^2-\lambda_{1,0}\norm{\varphi}_{-*}$, which tells the falseness of the expansion when $\varphi$ is not in $L^2\subset L^p$). We borrow a lemma due to Christ \cite{chr14} to get the lower bound, which Christ used to consider the stability for the Hausdorff-Young inequality functional.

\blem\label{christ}\emph{(Christ's lemma)}
For any exponents $p<2\le q$, there exist positive constants, $c_0, c_1, \eta_0$, $\gamma(>1)$, s.t., given any small $\eta,\delta$ satisfying $0< \eta \le \eta_0, ~0<\delta<\eta^\gamma$, any $(L^p,L^q)$-bounded linear operator for an arbitrary pair of measure spaces with norm $\norm{T}$ and some nonvanishing extremizer $F$ (i.e. $|TF|_q = \norm{T}|F|_p$), and any $f\in L^p$ satisfying $|f|_p\le \delta|F|_p, ~\re \int f\bar{F}|F|^{p-2}=0$, if we decompose $f=f_1+f_2$, where $f_1=f\chi_{|f|\le \eta}$, then we have the following estimate of the functional
\beq\label{christ-}\frac{|T(F+f)|_q}{\norm{T}|F+f|_p}\le 1 + \phi(f_1) +c_1\eta|f_1|_p^2|F|_p^{-2}-c_0\eta^{2-p}|f_2|_p^p|F|_p^{-p},\eeq where
\begin{align}
\phi(f)=\frac{q-1}{2}|TF|_q^{-q}\int\bra{\re\frac{Tf}{TF}}^2|TF|^q+\frac{1}{2}|TF|_q^{-q}\int\bra{\im\frac{Tf}{TF}}^2|TF|^q\nonumber\\
-\frac{p-1}{2}|F|_p^{-p}\int\bra{\re\frac{f}{F}}^2|F|^p-\frac{1}{2}|F|_p^{-p}\int\bra{\im\frac{f}{F}}^2|F|^p. \label{christ--}
\end{align}
\elem

Then we have the following local estimate:
\bpro\label{local-}\emph{(Local stability for HLS)}
Let $Q,s,p$ be the same as that in Theorem \ref{t-hls}, then there exist positive constants $\alpha_0,\alpha_1$ depending only on dimension $Q$ and $s$, s.t., $~\forall~ 0\not\equiv f\in L^p$,
\begin{align}\label{f-local-}
o\bra{\frac{d_p(f,\ca{M}_{-*})}{|f|_p}} + \alpha_1~ \frac{d_p(f,\ca{M}_{-*})}{|f|_p} \ge C^{-1}-\frac{\norm{f}_{-*}^2}{|f|_p^2} \ge \alpha_0~ \frac{d_p^2(f,\ca{M}_{-*})}{|f|_p^2} + o\bra{\frac{d_p^2(f,\ca{M}_{-*})}{|f|_p^2}}.
\end{align}
\epro
\bpf
We denote $f_p$ to be the $L^p$-nearest point, and locally, we may assume $f_p=c|J_\tau|^{1/p}$ for some $c\neq 0$ and conformal transformation $\tau$, noting $\frac{c}{|f|_p}\rightarrow 1$ as $\frac{d_p(f,\ca{M}_{-*})}{|f|_p}\rightarrow 0$. Denote $g=c^{-1}f\circ \tau^{-1} |J_{\tau^{-1}}|^{1/p}$, then $g_p=1$ and we write $g=1+\varphi$ with $|\varphi|_p\rightarrow 0$ and $\int\varphi=0$ (this condition need to be satisfied in Christ's lemma, and if $\int \varphi\neq0$, then update $\varphi$ by $\frac{\varphi-\dashint \varphi}{1+\dashint \varphi}$ and note that $|\varphi-\dashint \varphi|_p=d_p(g,1+\dashint \varphi)\ge d_p(g,\m_{-*})=|\varphi|_p$). Then it suffices to prove (\ref{f-local-}) for $g$ as the inequality is invariant under above transformation.
For lower bound, we use Lemma \ref{christ} (taking there $T=\ca{A}_s^{-1/2}, q=2, p=\frac{2Q}{Q+s}, F=1, \norm{T}=C^{-1/2}, f=\varphi, \delta=|\varphi|_p/|1|_p, \eta=(1+\epsilon)\delta^{1/\gamma}$), then from (\ref{christ-},\ref{christ--}) we have
\[\frac{C\norm{g}_{-*}^2}{|g|_p^2} \le \bra{1+\frac{1}{2}\norm{1}_{-*}^{-2}\norm{\varphi_1}_{-*}^2-\frac{p-1}{2}|1|_p^{-p}|\varphi_1|_2^2 +c_1\eta|\varphi_1|_p^2|1|_p^{-2}-c_0\eta^{2-p}|\varphi_2|_p^p|1|_p^{-p}}^2.\]
We easily see $\norm{\varphi_1}_{-*}^2\lesssim |\varphi_1|_p^2 = o(|\varphi|_p^2), |\varphi_1|_2^2\le \eta^{2-p}|\varphi|_p^p$ and it follows from $\eta\gtrsim |\varphi|_p$ that
\[1-\frac{C\norm{g}_{-*}^2}{|g|_p^2}\gtrsim \eta^{2-p}|\varphi|_p^p+o(|\varphi|_p^2) \gtrsim |\varphi|_p^2+o(|\varphi|_p^2),\] which gives the right side of (\ref{f-local-}).
For the upper bound, we denote the $H^{-*}$-nearest point $g_*=c'|J_\tau'|^{1/p}$, where $c'\neq 0$ as $g_*\xrightarrow{\norm{\cdot}_{-*}} 1$. Then from the conformal invariance and extremizer property, we have \[|g|_p^2-C\norm{g}_{-*}^2 \le |g|_p^2-C\norm{g_*}_{-*}^2=|g|_p^2-|g_*|_p^2 = |g|_p^2-|1|_p^2+(1-|c'|^2)|1|_p^2.\] Note that, $|g|_p^2-|1|_p^2\le |\varphi|_p(|g|_p+|1|_p) =2|1|_p|\varphi|_p+o(|\varphi|_p)$ and $|c'|\rightarrow 1$, actually, $\abs{1-|c'|}\lesssim |1-g_*|_{-*}\lesssim |\varphi|_p+|g-g_*|_{-*} \lesssim |\varphi|_p$. Above all, we get
\[|g|_p^2-C\norm{g}_{-*}^2 \lesssim |\varphi|_p +o(|\varphi|_p),\] then we get a local upper bound for $f$ by substituing back,
\[1- \frac{C\norm{f}_{-*}^2}{|f|_p^2} \lesssim \frac{c}{|f|_p} \bra{\frac{d_p(f,\m_{-*})}{|f|_p}+o\bra{\frac{d_p(f,\m_{-*})}{|f|_p}}}\] and finally from $\frac{|c|}{|f|_p}\lesssim 1+\abs{1-\frac{|c|}{|f|_p} |1|_p} \lesssim 1+ \abs{\frac{f}{|f|_p}-\frac{f_p}{|f|_p}}_p = 1+ \frac{d_p(f,\m_{-*})}{|f|_p}$, we get the left side of (\ref{f-local-}) (we can also see directly from the conformal invariance).  Proposition (\ref{local-}) is then proved.
\epf

For the global stability, we only need the recovery of compactness (dual of Lemma \ref{compact}, see a different formulation in Lemma 4.6 of \cite{fl12} and we can also use suitably adapted concentration compactness argument) and then apply the contradiction argument.

\blem\label{compact-}\emph{(Recovery of compactness for HLS)}
Let $(f_j)$ be an (non-vanishing) extremizing sequence of inequality \emph{(\ref{hls})} (or its functional) on the sphere, i.e., $\frac{\norm{f_j}_{-*}^2}{|f_j|_p^2}\xrightarrow{j\rightarrow\infty} C^{-1}$, then  $\frac{d_p(f_j,\ca{M}_{-*})}{|f_j|_p} \xrightarrow{j\rightarrow \infty} 0$.
\elem
\emph{Proof of Theorem \ref{t-hls}.}
We assume the statement for formula (\ref{f-hls-}) in Theorem \ref{t-hls} is wrong, then there exist two non-degenarate sequences of functions $(f_j)$ and $(g_j)$ satisfying
\beq\label{contradiction} \lim_{j\rightarrow \infty}\frac{C^{-1}-\frac{\norm{f_j}_{-*}^2}{|f_j|_p^2}}{\frac{d_p^2(f_j,\m_{-*})}{|f_j|_p^2}} = 0, \quad~\text{and}~\quad \lim_{j\rightarrow \infty}\frac{C^{-1}-\frac{\norm{g_j}_{-*}^2}{|g_j|_p^2}}{\frac{d_p(g_j,\ca{M}_{-*})}{|g_j|_p}} = \infty.\eeq
Then $\lim_{j\rightarrow \infty} \frac{\norm{f_j}_{-*}^2}{|f_j|_p^2}= C^{-1}$, recalling $\frac{d_p(f,\ca{M}_{-*})}{|f|_p}\le 1$, which gives $\lim_{j\rightarrow \infty} \frac{d_p(f_j,\ca{M}_{-*})}{|f_j|_p}=0$ from Lemma \ref{compact-}. Simultaneously, we have $\lim_{j\rightarrow \infty} \frac{d_p(g_j,\ca{M}_{-*})}{|g_j|_p}=0$. Actually, we see this from Lemma \ref{compact-} if $\lim_{j\rightarrow \infty} \frac{\norm{g_j}_{-*}^2}{|g_j|_p^2}= C^{-1}$, and otherwise after moving to a subsequence as $\liminf_{j\rightarrow \infty} \frac{\norm{g_j}_{-*}^2}{|g_j|_p^2}< C^{-1}$. So, the formulas in (\ref{contradiction}) are local limitations, which contradicts with the local bound (\ref{f-local-}) in Proposition \ref{local-}. \qed\\

\subsection{Proof of Proportional Inequalities Between Dual Remainder Terms}

Now we are going to prove the relations between remainder terms of FS and HLS inequalities and that of BO and Log-HLS inequalities. First, for the FS and HLS.

\emph{Proof of Theorem \ref{t-relation}.}
For the local bound (\ref{local bound}), we expand the two remainder terms by Taylor expansion. Take $f_*$ to be the $H^*$-nearest point, and as before we may assume $f_*\not\equiv 0$ (we can add other condition like $\norm{f}_*=1$ locally), and further $f_*=1, f=1+\varphi$, where $\varphi\perp T_1\m_*$ under $\ip{\cdot,\cdot}_*$. We denote the two remainder functionals respectively by
\[I_1(f)=\norm{f}_*^2-C|f|_q^2, \quad I_2(f)=|f|_p^2-C\norm{f}_{-*}^2.\]  Then we expand and estimate the two terms
\begin{align*}
I_1(f)=&~ I_1(1+\varphi)\\
=&~I_1(1)+\frac{d}{dt}\Big|_{t=0}I_1(1+t\varphi)+\frac{1}{2}\frac{d^2}{dt^2}\Big|_{t=0}I_1(1+t\varphi)+o(\norm{\varphi}_*^2)\\
=&~  \norm{\varphi}_*^2-(q-1)C|\s|^{2/q-1}|\varphi|_2^2+o(\norm{\varphi}_*^2)\\
I_2(f^{q/p})=&~I_1((1+\varphi)^{q/p})\\
=&~I_2(1)+\frac{d}{dt}\Big|_{t=0}I_2((1+t\varphi)^{q/p})+\frac{1}{2}\frac{d^2}{dt^2}\Big|_{t=0}I_2((1+t\varphi)^{q/p})+o(\norm{\varphi}_*^2)\\
=&~(q-1)^2|\s|^{2/p-1}|\varphi|_2^2-\frac{q}{p}\bra{\frac{q}{p}-1}C<1,\varphi^2>_{-*}-\bra{\frac{q}{p}}^2C\norm{\varphi}_{-*}^2+o(\norm{\varphi}_*^2)\\
=&~\bra{(q-1)^2|\s|^{2/p-1}-\frac{q}{p}\bra{\frac{q}{p}-1}C\lambda_{0,0}^{-1}}|\varphi|_2^2-\bra{\frac{q}{p}}^2C\norm{\varphi}_{-*}^2+o(\norm{\varphi}_*^2)\\
=&~(q-1)|\s|^{2/p-1}|\varphi|_2^2-\bra{\frac{q}{p}}^2C\norm{\varphi}_{-*}^2+o(\norm{\varphi}_*^2),
\end{align*}
where we have used $C\lambda_{0,0}^{-1}=|\s|^{2/p-1}$ and the $C^2(H^*)$ property of $I_1(f)$ and $I_2(f^{q/p})$ (this can be checked easily as in the proof of (\ref{taylor}), Proposition \ref{local}). Take bispherical harmonics decomposition on $\varphi$ and recalling $\varphi\perp \ca{H}_0\oplus\ca{H}_1$, we can set $\varphi=\sum_{j+k\ge 2} \varphi_{j,k}$, then we want to study the following quotient:
\begin{align*}
\frac{I_2(f^{q/p})}{I_1(f)}=(q-1)|\s|^{2/p-1}
\frac{\sum_{j+k\ge 2} \bra{1-\frac{q}{p}C|\s|^{1-2/p}\lambda_{j,k}^{-1}}|\varphi_{j,k}|_2^2}
{\sum_{j+k\ge 2} \bra{\lambda_{j,k}-(q-1)C|\s|^{2/q-1}}|\varphi_{j,k}|_2^2}+o(1)
\end{align*}
Note that \[\frac{1-\frac{q}{p}C|\s|^{1-2/p}\lambda_{j,k}^{-1}}{\lambda_{j,k}-(q-1)C|\s|^{2/q-1}}=\frac{1}{\lambda_{j,k}},\] we then have
\[\limsup_{0<d(f,\m_*)\rightarrow 0, d(f,\m_*)<\norm{f}_*}\frac{I_2(f^{q/p})}{|f|_q^{2(q-2)}I_1(f)}= \frac{(q-1)|\s|^{1-2/p}}{\lambda_{2,0}}=\frac{\lambda_{1,0}}{C\lambda_{2,0}}=\frac{1}{C}\frac{Q+4-s}{Q+4+s}.\]
For the global bound (\ref{global bound}), we use the square idea as in \cite{dj14}, for $r>0$ to be fixed later,
\begin{align*}
0&\le\abs{|f|_q^r\ca{A}_s^{1/2}f-C\ca{A}_s^{-1/2}f^{q/p}}_2^2\\
&=|f|_q^{2r}\abs{\ca{A}_s^{1/2}f}_2^2+C^2\abs{\ca{A}_s^{-1/2}f^{q/p}}_2^2-2C|f|_q^r\ip{\ca{A}_s^{1/2}f,\ca{A}_s^{-1/2}f^{q/p}}\\
&=|f|_q^{2r}\bra{\norm{f}_*^2-C|f|_q^{-r+q}}-C\bra{\abs{f^{q/p}}_p^{p(r+q)/q}-C\norm{f^{q/p}}_{-*}^2}
\end{align*}
Take $r=q-2=\frac{2s}{Q-s}$, then $-r+q=p(r+q)/q=2$, and
\[|f|_q^{2r}\bra{\norm{f}_*^2-C|f|_q^2}\ge C\bra{\abs{f^{q/p}}_p^{2}-C\norm{f^{q/p}}_{-*}^2}.\] So Theorem \ref{t-relation} is proved. \qed\\

Now we are going to prove the endpoint case --- relation between remainder terms of BO and Log-HLS inequalities, using the result about FS and HLS by differentiation argument globally and still the Taylor expansion for the local estimate.

\emph{Proof of Theorem \ref{t-relation-}.}
For global bound (\ref{global-limit}), we can see from the remainder terms control between FS and HLS inequality (\ref{global bound})  in Theorem \ref{t-relation}) by the differential argument. Actually we can write the inequality (\ref{global bound})  to the following form
\[I_1(f)\ge I_2(f),\]
where the associated two functionals are
\begin{align*}
I_1(f)=&~ \bra{\dashint |f|^q}^{2s/Q}\bra{\dashint f\ca{A}_sf-\frac{\Gamma(\frac{Q+s}{4})}{\Gamma(\frac{Q-s}{4})}\bra{\dashint |f|^q}^{2/q}},\\
I_2(f)=&~ \frac{\Gamma^2(\frac{Q+s}{4})}{\Gamma^2(\frac{Q-s}{4})}\bra{\bra{\dashint |f|^q}^{2/p}-\frac{\Gamma^2(\frac{Q+s}{4})}{\Gamma^2(\frac{Q-s}{4})} \dashint f^{q/p}\ca{A}_s^{-1}f^{q/p}},
\end{align*}
then from (\ref{dif}), $\ca{A}_s=\frac{\Gamma(\frac{Q+s}{4})}{\Gamma(\frac{Q-s}{4})}\ca{A}_Q'+o(\lambda)$ and by Taylor expansion on $\lambda$, we get
\begin{align*}
I_1\left(1+\frac{\lambda}{2Q}f\right)=&~ \bra{\bra{\dashint e^f}^2+o(1)}\Bigg(\dashint 1\ca{A}_s1+\bra{\frac{\lambda}{2Q}}^2\dashint \Big(f-\dashint f\Big)\ca{A}_s\Big(f-\dashint f\Big)\\
&~\qquad~ \qquad~ \qquad -\frac{\Gamma^2(\frac{Q+s}{4})}{\Gamma^2(\frac{Q-s}{4})} \bra{1+\frac{\lambda}{Q}\log \dashint e^{f-\dashint f}+o(\lambda)}\Bigg)\\
=&~ \frac{1}{Q}\bra{\frac{n!}{4}}^2  \bra{\dashint e^f}^2\bra{\frac{1}{2(n+1)!}\dashint f\ca{A}_Q'f-\log \dashint e^{f-\dashint f}}\lambda^3+o(\lambda^3),
\end{align*} and
\begin{align*}
I_2\bra{1+\frac{\lambda}{2Q}f}=&~ \bra{\frac{n!\lambda}{4}}^2\Bigg(\bra{\dashint e^f}^2\bra{1-\Big(\frac{\log\dashint e^f}{Q}+\frac{1}{2Q}\frac{\dashint e^f f^2}{\dashint e^f}\Big)\lambda+o(\lambda)}\\
&~  - \bra{\dashint \bra{1+\frac{\lambda}{2Q}f}^{\frac{2Q-\lambda}{\lambda}}}^2- \frac{\Gamma(\frac{Q+s}{4})}{\Gamma(\frac{Q-s}{4})}\dashint\Big(e^f-\dashint e^f\Big)\ca{A}_Q'^{-1}\Big(e^f-\dashint e^f\Big)\Bigg)  \\
=&~ \frac{1}{Q}\bra{\frac{n!}{4}}^2  \Bigg(\dashint e^f\dashint e^ff-\bra{\dashint e^f}^2\log \dashint e^f - \frac{(n+1)!}{2}\\
&~ \qquad ~\qquad \times\dashint \Big(e^f-\dashint e^f\Big)\ca{A}_Q'^{-1}\Big(e^f-\dashint e^f\Big)\Bigg)\lambda^3+o(\lambda^3).
\end{align*} Then (\ref{global-limit}) is proved by comparing the dominating terms: dividing the two formulas by $\lambda^3$ and taking $\lambda=Q-s\rightarrow 0$ .\\
For local estimate (\ref{local-limit}), by constant translation invariance of the two terms in the quotient formula, we can assume $f\perp \ca{H}_0\oplus\ca{H}_1$. Now, we consider Taylor expansion of the quotient, which we denote by $I(f)$,
\begin{align*}\
I(f)=&~ \frac{\frac{1}{2(n+1)!}(\dashint f\ca{A}_Q'f-(n+1)!\dashint f^2)+o(\dashint f \ca{A}_Q'f)}{\frac{1}{2}(\dashint f^2-(n+1)!\dashint f\ca{A}_Q'^{-1}f)+o(\dashint f\ca{A}_Q'f)}\\
=&~\frac{1}{(n+1)!}\frac{\sum_{j\ge 2}(\frac{\Gamma(j+n+1)}{\Gamma(j)}-(n+1)!)(|f_{j,0}|_2^2+|f_{0,j}|_2^2)}{\sum_{j\ge 2}(1-\frac{\Gamma(j)}{\Gamma(j+n+1)}(n+1)!)(|f_{j,0}|_2^2+|f_{0,j}|_2^2)}\\
\ge &~ n+2.
\end{align*} \qed

\section*{Acknowledgement}
The work was partially done during An Zhang's visit in the department of mathematics, UC, Berkeley, as a visiting student researcher. He would like to thank the department for hospitality and especially professor Michael Christ for leading him to this interesting research field, for his host at UC, Berkeley, and discussions and valuable comments which improve the final version of this note.

\bibliographystyle{amsalpha}
\bibliography{reference}

\vskip 3\baselineskip
\flushleft

Heping Liu\\
School of Mathematical Science,
Peking University;
Beijing, China.\\
\emph{email:} \texttt{hpliu@math.pku.edu.cn}
\vskip 2\baselineskip
An Zhang\\
School of Mathematical Science,
Peking University;
Beijing, China.\\
\emph{email:} \texttt{anzhang@pku.edu.cn}\\
Current Address:\\
Department of Mathematics,
University of California, Berkeley;
Berkeley, CA, USA. \\
\emph{email:} \texttt{anzhang@math.berkeley.edu}

\end{document}